\documentclass[reqno,12pt]{amsart}
\usepackage{amsmath,amsfonts,mathrsfs}%amssymb,
\newtheorem{thm}{Theorem}
\newtheorem{cor}[thm]{Corollary}

\newtheorem{prop}[thm]{Proposition}

\newtheorem{lthm}{Theorem A\!\!}

\newtheorem{llem}[lthm]{Lemma A\!\!}
\newtheorem{lprop}[lthm]{Proposition A\!\!}

\newtheorem{defi}{Definition}

\def\mysection#1{\section{#1} \setcounter{equation}{0}}

\def\tr{\mathop{\rm tr}\nolimits}

\def\min{\mathop{\rm min}}
\def\max{\mathop{\rm max}}

\newcommand{\norm}[1]{\left\Vert#1\right\Vert}

\newcommand{\Real}{\mathbb R}
\newcommand{\NN}{\mathbb N}
\newcommand{\ZZ}{\mathbb Z}

\title[Gibbs-like measure for spectrum]
{Gibbs-like measure for spectrum of\\ a class of one-dimensional
Schr\"odinger operator\\ with Sturm potentials }
\author[FAN, LIU, WEN]{Shen FAN$^\dagger$, Qing-Hui LIU$^*$\and
Zhi-Ying WEN$^{\dagger\ddag}$}
\thanks{$^*$ Department of Computer Science and Engineering ,
Beijing Institute of Technology, Beijing, China. Supported by
National Natural Science Foundation of China No.10971013 and
Excellent young scholars Research Foundation of Beijing Institute of
Technology. }
\thanks{$^\dagger$ Department of Mathematics, Tsinghua University,
Beijing, 100084, P. R. China.}
\thanks{$^\ddag$
Supported by the National Nature Science Foundation of China
No.10631040 and National Basic Research Program of China(973
Program)2007CB814800.}
\date{}

\begin{document}

\begin{abstract}
Let $\alpha\in(0,1)$ be an irrational, and $[0;a_1,a_2,\cdots]$ the
continued fraction expansion of $\alpha$. Let $H_{\alpha,V}$ be the
one-dimensional Schr\"odinger operator with Sturm potential of frequency $\alpha$.
Suppose the potential strength $V$ is large enough and $(a_i)_{i\ge1}$ is bounded.
We prove that the spectral generating bands possess  properties of bounded distortion,
bounded covariation and there exists Gibbs-like measure on the spectrum $\sigma(H_{\alpha,V})$.
As an application, we prove that
$$\dim_H \sigma(H_{\alpha,V})=s_*,\quad \overline{\dim}_B \sigma(H_{\alpha,V})=s^*,$$
where $s_*$ and $s^*$ are lower and upper pre-dimensions.
\end{abstract}

\maketitle

\mysection{Introduction} Since the discovery of quasi-crystal by
Schechtman et al. (\cite{SBGC}), the one dimensional discrete
Schr\"odinger operators with Sturm potentials have largely been
studied, see \cite{BIST,R,Su} and references therein.
The discrete Schr\"odinger operator acting on $l^2(\mathbb{Z})$ is defined as
follows: for any $\psi=\{\psi_n\}_{n\in\ZZ}\in l^2(\mathbb{Z})$,
\begin{equation}\label{schr}
(H_{\alpha,V}\psi)_n=\psi_{n-1}+\psi_{n+1}+v_n
\psi_n,\ \forall
n\in\ZZ.
\end{equation}
The potential $(v_n)_{n\in\mathbb{Z}}$ we discuss in this paper is the Sturm potential, i.e.,
\begin{equation}\label{sturm}
v_n=V\chi_{[1-\alpha,1)}(n\alpha+\phi \mod 1),\quad \forall n\in\mathbb{Z},
\end{equation}
where $\alpha\in(0,1)$ is an irrational, and is called frequency,
$V>0$ is called potential strength or coupling, $\phi\in[0,1)$ is
called phase. We will study the structure of the spectrum of the
operator which we denote by $\sigma(H_{\alpha,V})$, in particular
the fractal dimensions of $\sigma(H_{\alpha,V})$. It is known that
$\sigma(H_{\alpha,V})$ is independent of phase $\phi$, we set
$\phi=0$.

It is proved by Bellissard, Iochum, Scoppola and Testart (\cite{BIST}, 1989) that
$\sigma(H_{\alpha,V})$ is a Cantor set of Lebesgue measure zero.
(On the other direction, in stead of Sturm potential, some authors
considered the primitive substitutive potential, it is proved that
in this case, the spectrum  is also a Cantor set of Lebesgue
measure zero. For more details, we refer to \cite{Lenz},\cite{LTWW}.)

Since then, whether the Hausdorff dimension of $\sigma(H_{\alpha,V})$
is strictly less than $1$ and strictly greater than $0$ have
absorbed a lot of attentions.
Raymond \cite{R}(1997) studied this
problem under the restriction $V>4$. Connected with the continued
fraction expansion of $\alpha$, he exhibited an interesting
recurrent structure of the spectrums. And for
$\alpha=\frac{\sqrt{5}-1}{2}$, i.e., the golden mean, he gave an
upper bound of the Hausdorff dimension of the corresponding
spectrum, which is strictly less than $1$.

Damanik, Killip and Lenz \cite{DKL}(2000) proved that if $\alpha$
has bounded density (this means if $[0;a_1,a_2,\cdots]$ is the
continued fraction expansion of $\alpha$, then
$\limsup\limits_{k\rightarrow\infty}\frac{1}{k}\sum_{i=1}^k
a_i<\infty$), then the Hausdorff dimension of the spectral measure
of $H_{\alpha,V}$ is strictly greater than $0$. Since the spectral
measure is supported by the spectrum $\sigma(H_{\alpha,V})$, the
Hausdorff dimension of the spectrum has also strictly positive
lower bound.

To estimate fractal dimensions of the spectrum of $H_{\alpha,V}$,
one of the key steps is to estimate the length of spectral
generating bands. Raymond \cite{R} has treated the case of frequency
$\alpha$ being golden mean with $V>4$. Based on the Raymond's
method, for all irrational frequency and $V>20$, Liu and Wen
(\cite{LW},\cite{LW05}) established multi-type Moran construction
among different spectral generating bands, developed a very fine
estimating technique for the length of the bands of different orders
of the spectrum $\sigma(H_{\alpha,V})$, and generalized some
techniques analogous to the studies of Moran structure in
\cite{FWW,HRWW}, they proved the following result.

\begin{lthm}\label{main}
Let $\alpha=[0;a_1,a_2,\cdots]$ be irrational and
$$K=\liminf_{k\rightarrow\infty}(a_1a_2\cdots a_k)^{1/k}.$$
Let $V>20$ and $t_1=\frac{3}{V-8}$, $t_{2}=\frac{1}{4(V+8)}$.

{\rm(1)}\ If $K<\infty$, then
$$\max\{\frac{\ln 2}{10\ln2-3\ln t_{2}},
\frac{\ln K-\ln 3}{\ln K-\ln (t_{2}/3)}\}\leq
\dim_H\,\sigma(H_{\alpha,V}) \leq \frac{2\ln K+\ln3}{2\ln K-\ln t_1};$$

{\rm(2)}\ If $K=\infty$, then
$$\dim_H\,\sigma(H_{\alpha,V})=1.$$
\end{lthm}

Note that this theorem implies that if $K<\infty$, then
$\dim_H(\sigma(H_{\alpha,V}))$ tends to 0 when $V$ tends to
infinity.

Damanik, Embree, Gorodetski, and Tcheremchantsev (\cite{DEGT})
proved that, for golden mean $\alpha$,
\begin{equation*}
\lim_{V\to \infty} (\log V)\, \overline{\dim}_B \sigma(H_{\alpha,V})
= \log (1+\sqrt{2}),
\end{equation*}
and found that $\dim_H\, \sigma(H_{\alpha,V})=\overline{\dim}_B\, \sigma(H_{\alpha,V})$
by applying dimensional theory of dynamical system.

Liu, Peyri\`ere, Wen\cite{LPW07} extended their results to case of
$\alpha=[0;a_1,a_2,\cdots]$ with $(a_n)_{n\ge1}$ bounded. They
proved that, for pre-dimensions $0\le s_*\le s^*\le1$ (which will be defined
in \S2),
$${\dim}_H\, \sigma(H_{\alpha,V})\le s_*,\quad \overline{\dim}_B \sigma(H_{\alpha,V})
\ge s^*,$$ and
\begin{equation*}
\lim_{V\to \infty} s_*\log V
= -\log f_*(\alpha),\quad
\lim_{V\to \infty} s^*\log V
= -\log f^*(\alpha),
\end{equation*}
where $f_*(\beta)$ and $f^*(\beta)$ are the positive roots of the
equations
\begin{equation*}\begin{array}{l}
\liminf_{n\rightarrow\infty} \norm{{\mathbf R}_1(x){\mathbf
R}_2(x)\cdots {\mathbf R}_n(x)}^{1/n}=1\\
\limsup_{n\rightarrow\infty} \norm{{\mathbf R}_1(x){\mathbf
R}_2(x)\cdots {\mathbf R}_n(x)}^{1/n}=1,
\end{array}
\end{equation*}
and for any $0<x\le1$ and $n\ge1$,
\begin{equation*}
{\mathbf R}_n(x) := \begin{pmatrix}
0&x^{(a_n-1)}&0\\
(a_n+1)x&0&a_n x\\
a_n x&0&(a_n-1)x
\end{pmatrix}.
\end{equation*}
If $\alpha=[0;1,1,\cdots]$, then
$f_*(\alpha)=f^*(\alpha)=(1+\sqrt{2})^{-1}$, which is the positive
root of
$$\det\left(\left[\begin{array}{ccc}
0&1&0\\ 2x&0&x\\ x&0&0 \end{array}\right]-I\right)=-1+2x+x^2=0.$$
They also show that there are frequencies $\alpha$ with $f_*(\alpha)<f^*(\alpha)$.

In this paper,  we will consider the general formula of the
dimensions of the spectrum for the case $(a_n)_{n\ge1}$ being
bounded. For this aim, we establish first the properties of bounded
variation, bounded covariation for spectral generating bands, then
prove the existence of Gibbs-like measures for spectrum, finally we
give a general result of the Hausdorff dimension and upper box
dimension of the spectrum, that is,
$${\dim}_H\, \sigma(H_{\alpha,V})=s_*,\quad
\overline{\dim}_B\, \sigma(H_{\alpha,V})=s^*.$$

The remainder of the paper will be organized as follows: in Section
2, we introduce spectral structure and state the main results of the
paper; Section 3 will be devoted to the proofs of these results.

\mysection{Spectral structure}

We discuss first some facts on the structure of
$\sigma(H_{\alpha,V})$.

Let $\alpha=[a_1,a_2,\cdots,a_i,\cdots]\in(0,1)$ be an irrational,
let $p_k/q_k$$(k>0)$ be the $k$-th asymptotic fraction of $\alpha$
given by:
$$\begin{array}{l}
p_{-1}=1,\quad p_0=0,\quad p_{k+1}=a_{k+1} p_k+p_{k-1},\ k\ge 0,\\
q_{-1}=0,\quad q_0=1,\quad q_{k+1}=a_{k+1} q_k+q_{k-1},\ k\ge 0.
\end{array}$$

\smallskip

\noindent Let $k\geq1$ and $x\in\Real$, the transfer matrix $M_k(x)$
over $q_k$ sites is defined by
$${\mathbf M}_k(x):=
\left[\begin{array}{cc}x-v_{q_k}&-1\\ 1&0\end{array}\right]
\left[\begin{array}{cc}x-v_{q_k-1}&-1\\ 1&0\end{array}\right]
\cdots \left[\begin{array}{cc}x-v_2&-1\\ 1&0\end{array}\right]
\left[\begin{array}{cc}x-v_1&-1\\ 1&0\end{array}\right],$$
where $v_n$ is defined in \eqref{sturm} and by convention, take
$$\begin{array}{l}
{\mathbf M}_{-1}(x)= \left[\begin{array}{cc}1&-V\\
0&1\end{array}\right],\quad {\mathbf M}_{0}(x)=
\left[\begin{array}{cc}x&-1\\ 1&0\end{array}\right].
\end{array}$$

\smallskip

For $k\ge0$, $p\ge-1$, let $t_{(k,p)}(x)=\tr {\mathbf M}_{k-1}(x) {\mathbf M}_k^p(x)$ and
$\sigma_{(k,p)}=\{x\in\Real:|t_{(k,p)}(x)|\leq2\}$, where  $tr M$
stands for the trace of the matrix $M$.

With these notations, we
collect some known facts that will be used later, for more
details, we refer to \cite{BIST,R,Su,T}.
\begin{itemize}
\item[(A)]\ Renormalization relation.
For any $k\ge0$
\begin{equation}\label{renorm}
{\mathbf M}_{k+1}(x)={\mathbf M}_{k-1}(x)({\mathbf
M}_k(x))^{a_{k+1}},
\end{equation}
so, $t_{(k+1,0)}=t_{(k,a_{k})}$, $t_{(k,-1)}=t_{(k-1,a_k-1)}$.

\item[(B)]\ Structure of $\sigma_{(k,p)}(k\ge0,p\ge-1)$.

For $V>0$, $\sigma_{(k,p)}$ is made out of $\deg t_{(k,p)}$ separated
closed intervals.

\item[(C)]\ Trace relation.

By defining $\Lambda(x,y,z)=x^2+y^2+z^2-xyz-4$,
\begin{equation}\label{invariant}
\Lambda(t_{(k+1,0)},t_{(k,p)},t_{(k,p+1)})=V^2.
\end{equation}
Thus for any $k\in\NN$, $p\geq 0$ and $V>4$,
\begin{equation}\label{empty}
\sigma_{(k+1,0)}\cap \sigma_{(k,p)}\cap\sigma_{(k,p-1)}=\emptyset.
\end{equation}

\item[(D)]\ Covering property.

For any $k\ge0$, $p\ge-1$,
$$\sigma_{(k,p+1)}\subset \sigma_{(k+1,0)}\cup \sigma_{(k,p)},$$
then
$$(\sigma_{(k+2,0)}\cup\sigma_{(k+1,0)})\subset
(\sigma_{(k+1,0)}\cup\sigma_{(k,0)}).$$
Moreover
$$\sigma(H_{\alpha,V})=\bigcap_{k\ge0}(\sigma_{(k+1,0)}\cup\sigma_{(k,0)}).$$
\end{itemize}

We call the constructive intervals of $\sigma_{(k,p)}$ the {\em
bands}. When we discuss only one of these bands, we often denote it
as $B_{(k,p)}$. Property (B) also implies $t_{(k,p)}(x)$ is monotone
on $B_{(k,p)}$, and
$$t_{(k,p)}(B_{(k,p)})=[-2,2],$$
we call $t_{(k,p)}$ the {\em generating polynomial} of $B_{(k,p)}$.

\begin{defi}{\rm (\cite{R,LW})}
For $V>4$, $k\ge0$, we define three types of bands as follows:

$(k,{\rm I})$-type band: a band of $\sigma_{(k,1)}$ contained in a
band of $\sigma_{(k,0)}$;

$(k,{\rm II})$-type band: a band of $\sigma_{(k+1,0)}$ contained
in a band of $\sigma_{(k,-1)}$;

$(k,{\rm III})$-type band: a band of $\sigma_{(k+1,0)}$ contained
in a band of $\sigma_{(k,0)}$.
\end{defi}

The three kinds of types of bands are well defined(\cite{R}), and we
call these bands {\em spectral generating bands of order $k$} (the
type I band is called the type I gap in \cite{R}). Note that for
order $0$, there is only one $(0,{\rm I})$-type band
$\sigma_{(0,1)}=[V-2,V+2]$ (the corresponding generating polynomial
is $t_{(0,1)}=x-V$), and only one $(0,{\rm III})$ type band
$\sigma_{(1,0)}=[-2,2]$ (the corresponding generating polynomial is
$t_{(1,0)}=x$). They are contained in
$\sigma_{(0,0)}=(-\infty,+\infty)$ with corresponding generating
polynomial $t_{(0,0)}\equiv2$. For the convenience, we call
$\sigma_{(0,0)}$ the spectral generating band of order $-1$.

\smallskip

For any $k\ge-1$, denote by $\mathscr{G}_k$ the set of all spectral
generating bands of order $k$. By the properties (A),(B),(C) and (D),
for any $k\ge0$, we have
\begin{itemize}
\item $(\sigma_{(k+2,0)}\cup\sigma_{(k+1,0)})\subset
\bigcup_{B\in\mathscr{G}_k}B
\subset (\sigma_{(k+1,0)}\cup\sigma_{(k,0)})$,
and then
$$\sigma(H_{\alpha,V})=\bigcap_{k\ge-1}
\bigcup_{B\in\mathscr{G}_k}B;$$
\item
any $(k+1,I)$ or $(k+1,III)$-type band is contained in a $(k,II)$ or
$(k,III)$-type band; any $(k+1,II)$-type band is contained in a
$(k,I)$-type band;
\item
any $(k,II)$ or $(k,III)$ do not contain any $(k+1,II)$-type band;
any $(k,I)$-type band contain neither $(k+1,I)$ nor $(k+1,III)$-type
band.
\end{itemize}

To show that one band of order $k$ contains how many bands of order
$k+1$, we introduce Chebischev polynomial $S_p(x)$, which is defined
by
$$\begin{array}{l}
S_0(x)\equiv0,\quad S_1(x)\equiv1,\\
S_{p+1}(x)=x S_p(x)-S_{p-1}(x),\quad p\geq1.
\end{array}$$

By induction we see that
\begin{equation}\label{sin}
S_p(2\cos\theta)=\frac{\sin p\theta}{\sin\theta},\quad
\theta\in[0,\pi].
\end{equation}

Our study focus on the following three formulas according to the
types of the band(see \cite{BIST,R,LW}):
\begin{equation}\label{cII}
t_{(k,p)}=t_{(k,0)}S_{p+1}(t_{(k+1,0)})-t_{(k,-1)}S_p(t_{(k+1,0)}).
\end{equation}

\begin{equation}\label{cIII}
t_{(k,p+1)}=t_{(k,1)}S_{p+1}(t_{(k+1,0)})-t_{(k,0)}S_p(t_{(k+1,0)}).
\end{equation}

\begin{equation}\label{cI}
t_{(k,p+1)}=t_{(k+1,0)} t_{(k,p)}-t_{(k,p-1)}.
\end{equation}

These three formulas can be obtained by the following way: let $A$
be a $2\times2$ matrix with $|A|=1$, then by Caylay-Hamilton Theorem
$A^2-(\tr A) A+I=0$, and hence, for any $n>1$,
$$\begin{array}{rcl}
A^n&=&S_{n}(\tr A)\, A-S_{n-1}(\tr A)\, I\\
&=&S_{n+1}(\tr A)\, I-S_{n}(\tr A)\, A^{-1}.
\end{array}$$
Then take the trace in the both sides and by the definitions of
$t_{(k,p)}$, the three formulas come.

Now consider the equation
$$\Lambda(x,y,z)=V^2,$$
then
\begin{equation}\label{inv}
z_{\pm}(x,y,V)=\frac{xy}{2}\pm\frac{1}{2}\sqrt{4V^2+(4-x^2)(4-y^2)}.
\end{equation}
For two branches $z=z_+$ or $z=z_-$, let
$$\begin{array}{l}
z_1(x,y,V):=\frac{\partial z(x,y,V)}{\partial x},
\ z_2(x,y,V):=\frac{\partial z(x,y,V)}{\partial y},
\ z_{11}(x,y,V):=\frac{\partial^2 z(x,y,V)}{\partial x\partial x},\\[5pt]
\ z_{12}(x,y,V):=\frac{\partial^2 z(x,y,V)}{\partial x\partial y},
\ z_{21}(x,y,V):=\frac{\partial^2 z(x,y,V)}{\partial y\partial x},
\ z_{22}(x,y,V):=\frac{\partial^2 z(x,y,V)}{\partial y\partial y}.
\end{array}$$
For any $|x|\leq2$, $|y|\leq2$, and $V>4$,
\begin{equation}\label{dec}
\begin{array}{l}
V-2\le|z_\pm(x,y,V)|\le V+2,\quad|z_1(x,y,V)|\leq1,\ |z_2(x,y,V)|\leq1,\\
|z_{11}(x,y,V)|\leq1,\ |z_{12}(x,y,V)|\leq1,\ |z_{21}(x,y,V)|\leq1,\ |z_{22}(x,y,V)|\leq1.
\end{array}
\end{equation}

\medskip

In the previous papers \cite{R,LW}, the authors have estimated the
derivatives of the generating polynomials and the number of the
bands of different types through the formulas \eqref{cII},
\eqref{cIII} and \eqref{cI}. Since the present situation are much
more complicated, for treating the relations among different types
of bands, we introduce the notion of ladder as follows.

For any $n>k\ge-1$, let
$$B_n\subseteq B_{n-1}\subseteq\cdots\subseteq B_{k},$$
be a sequence of spectral generating bands from order $n$ to $k$. We
call the sequence $(B_i)_{i=k}^n$ an {\em initial ladder}, and the
bands $B_i(k\le i\le n)$ are called initial rungs. Now we are going
to modify the initial ladder by the following way: for any $i(k\le
i\le n-1)$,
\begin{itemize}
\item if $B_i$ is of $(i,I)$-type with $a_{i+1}=1$, we delete the rung
$B_{i+1}$(in this case $B_{i+1}$ must be $(i+1,II)$-type, then
$t_{(i+2,0)}=t_{(i,1)}$ and $t_{(i+1,-1)}=t_{(i,0)}$ implies
$B_{i+1}=B_i$);

\item if $B_i$ is of $(i,I)$-type with $a_{i+1}=2$, we change
nothing;

\item if $B_i$ is of $(i,I)$-type with $a_{i+1}>2$, we add rungs
$(B_{(i,p)})_{p=2}^{a_{i+1}-1}$ between $B_i$ and $B_{i+1}$ :
$$B_{i+1}=B_{(i,a_i)}\subset B_{(i,a_i-1)}\subset\cdots\subset
B_{(i,2)}\subset B_{(i,1)}=B_i;$$
\item if $B_i$ is of $(i,II)$ or $(i,III)$-type, we change nothing.
\end{itemize}
We get by this way a unique modified ladder which we relabel as
$$B_n=\hat{B}_m\subset\cdots\subset\hat{B}_1\subset \hat{B}_{0}=B_{k}.$$
We call $(\hat{B}_i)_{i=0}^m$ the {\em modified ladder}, and we
denote the corresponding generating polynomials by
$(\hat{h}_i)_{i=0}^m$. Note that any two consecutive initial rungs
can not be of type $I$ simultaneously,  so the length of the
modified ladder is larger than $[(n-k)/2]$.

Although we do not define type for the bands of order $-1$, note
that the bands of order $0$ are either of type $I$ or of type $III$,
we can view $B_{-1}=\sigma_{(0,0)}$ as a band of type $II$ or $III$, thus
we need not add rungs between $B_0$ and $B_{-1}$.

%Let
%$$B_n\subset B_{n-1}\subset\cdots\subset B_0\subset B_{-1}$$
%be spectral generating bands from order $n$ to order $-1$. If
%$n<m<k<l\le-1$, the overlap of two ladders $(B_i)_{i=k}^n$ and
%$(B_i)_{i=l}^m$ are coincide, so they can be conjuncted, the
%conjuncted ladder is just the ladder for $(B_i)_{i=l}^n$.

\medskip

Let $(\hat{B}_i)_{i=0}^m$ be a modified ladder and
$(\hat{h}_i)_{i=0}^m$ the correspondent generating polynomials.
Using the notion of ladders, we can unify the formulas
\eqref{cII}--\eqref{cI} by one single formula. To see this, for any
$-1<i<m$, note that $\hat{B}_i$ may be in one of the following four
situation: type $I$, $II$, $III$ of a order $k\ge0$ and an added
rung of a order $k\ge0$, and we distinguish them further into the
following three cases.

Case 1)\ \ $\hat{B}_i$ is of $(k,II)$-type.

In this case, $\hat{h}_i=t_{(k+1,0)}$, $\hat{h}_{i+1}=t_{(k,p)}$ for some $p\ge1$,
$\hat{h}_{i-1}=t_{(k,-1)}=t_{(k-1,a_k-1)}$ (note that
$\hat{B}_{i-1}$ is an added rung if $a_{k}>2$, this is also an
advantage to apply modified ladder ). We have
\begin{equation*}\label{kii}
\begin{array}{rcl}
\hat{h}_{i+1}=t_{(k,p)}&=&t_{(k,0)}S_{p+1}(t_{(k+1,0)})-t_{(k,-1)}S_p(t_{(k+1,0)})\\
&=&t_{(k,0)}S_{p+1}(\hat{h}_{i})-\hat{h}_{i-1}S_p(\hat{h}_i)\\
t_{(k,0)}&=&z_{\pm}(\hat{h}_{i},\hat{h}_{i-1},V).
\end{array}\end{equation*}

Case 2)\ \ $\hat{B}_i$ is of $(k,III)$-type.

In this case, $\hat{h}_i=t_{(k+1,0)}$, $\hat{h}_{i-1}=t_{(k,0)}$,
$\hat{h}_{i+1}=t_{(k,p+1)}$ for some $p\ge0$. We have
\begin{equation*}\label{kiii}
\begin{array}{rcl}
\hat{h}_{i+1}=t_{(k,p+1)}&=&t_{(k,1)}S_{p+1}(t_{(k+1,0)})-t_{(k,0)}S_p(t_{(k+1,0)})\\
&=&t_{(k,1)}S_{p+1}(\hat{h}_{i})-\hat{h}_{i-1}S_p(\hat{h}_i)\\
t_{(k,1)}&=&z_{\pm}(\hat{h}_{i},\hat{h}_{i-1},V),
\end{array}\end{equation*}

Case 3)\ \ $\hat{B}_i$ is of $(k,I)$-type, or,  an added rung in order $k$.

In this case, there exists $1\le p\le a_{k+1}$ such
that $\hat{h}_i=t_{(k,p)}$, $\hat{h}_{i+1}=t_{(k,p+1)}$,
$\hat{h}_{i-1}=t_{(k,p-1)}$, and
\begin{equation*}\label{ki}
\begin{array}{rcl}
\hat{h}_{i+1}=t_{(k,p+1)}&=&t_{(k+1,0)}t_{(k,p)}-t_{(k,p-1)}\\
&=&t_{(k+1,0)}S_2(\hat{h}_{i})-\hat{h}_{i-1}S_1(\hat{h}_i)\\
t_{(k+1,0)}&=&z_{\pm}(\hat{h}_i,\hat{h}_{i-1},V),
\end{array}\end{equation*}

We summarize the above three cases by
\begin{equation}\label{ladder-i}
\hat{h}_{i+1}(x)=z_{\pm}(\hat{h}_i(x),\hat{h}_{i-1}(x),V)
S_{p_i+1}(\hat{h}_{i}(x))-\hat{h}_{i-1}(x)S_{p_i}(\hat{h}_i(x)),
\end{equation}
where $p_i$ take values as follows,
\begin{equation}\label{pvalue}
p_i=\left\{
\begin{array}{cl}
a_{k+1},& \mbox{ if $\hat{B}_i$ is of $(k,III)$-type and $\hat{B}_{i+1}$ is of $(k+1,I)$-type,}\\
a_{k+1}-1,& \mbox{ if $\hat{B}_i$ is of $(k,III)$-type and $\hat{B}_{i+1}$ is of $(k+1,III)$-type,}\\
a_{k+1}+1,& \mbox{ if $\hat{B}_i$ is of $(k,II)$-type and $\hat{B}_{i+1}$ is of $(k+1,I)$-type,}\\
a_{k+1},& \mbox{ if $\hat{B}_i$ is of $(k,II)$-type and $\hat{B}_{i+1}$ is of $(k+1,III)$-type,}\\
1,& \mbox{ if $\hat{B}_i$ is of $(k,I)$-type or an added rung at order $k$.}
\end{array}
\right.
\end{equation}

\begin{defi}
For $p\ge1$, $1\le l\le p$, set
$$I_{p,l}=\left\{2\cos\frac{l+c}{p+1}\pi\ :\ |c|\le\frac{1}{10},
\ |S_{p+1}(2\cos\frac{l+c}{p+1}\pi)|\le\frac{1}{4}\right\}.$$
\end{defi}

By the definition, for any $1\le l\le p$, we have
$S_{p+1}(2\cos\frac{l\pi}{p+1})=0$;
$|S_{p+1}(2\cos\frac{l+c}{p+1}\pi)|\le\frac{1}{4}$ implies
$|c|\le\frac{1}{10}$; $\{I_{p,l}\}_{l=1}^{p}$ are $p$ disjoint
intervals in $[-2,2]$.

The following property comes from essentially \cite{R} and
\cite{LW}, for the completeness, we give a proof here.
\begin{lprop}\label{index}
Assume $V>20$. Let $(\hat{B}_i)_{i=0}^m$ be a modified ladder,
$(\hat{h}_i)_{i=0}^{m}$ the corresponding generating polynomials,
and $(p_i)_{i=1}^{m-1}$ be given as in \eqref{pvalue}. Then for any
$0<i<m$, there exist a unique $l(1\le l\le p)$ such that
$$\hat{h}_i(\hat{B}_{i+1})\subset I_{p_i,l}.$$
\end{lprop}
\begin{proof}
For convenience, we denote
$z_{\pm}(\hat{h}_i(x),\hat{h}_{i-1}(x),V)$ by $z_\pm(x)$.

Note that $S_{p}^2-1=S_{p-1}S_{p+1}$. For $\delta=\pm1$, by
\eqref{ladder-i} and a direct computation,
\begin{equation}\label{zero}
(S_p(\hat{h}_i)+\delta)(\hat{h}_{i+1}+\delta\hat{h}_{i-1})=
S_{p+1}(\hat{h}_i)\left(\left(z_\pm(x) S_p(\hat{h}_i)-
\hat{h}_{i-1}S_{p-1}(\hat{h}_i)\right)+\delta z_\pm(x)\right).
\end{equation}
Notice first for any $x\in\hat{B}_{i+1}$, we have $|z_\pm(x)|\ge
V-2$. On the other hand, it can be verified that
$$\Lambda(\hat{h}_{i+1},\hat{h_i},z_\pm(x) S_p(\hat{h}_i)
-\hat{h}_{i-1}S_{p-1}(\hat{h}_i))=V^2,$$ so for any
$x\in\hat{B}_{i+1}$ we also have
$$|z_\pm(x) S_p(\hat{h}_i)-\hat{h}_{i-1}S_{p-1}(\hat{h}_i)|\ge V-2.$$

Choosing suitably $\delta=1$ or $-1$ so that for any
$x\in\hat{B}_{i+1}$,
$$\left|\left(z_\pm S_p(\hat{h}_i)- \hat{h}_{i-1}S_{p-1}(\hat{h}_i)\right)+
\delta z_\pm\right|\ge2(V-2).$$

Since $\hat{h}_{i+1}$ and $\hat{h}_{i-1}$ are monotone on
$\hat{B}_{i+1}$, and
$$\hat{h}_{i+1}(\hat{B}_{i+1})=[-2,2],\quad \hat{h}_{i-1}(\hat{B}_{i+1})\subset[-2,2],$$
there exists a unique point $x_0\in\hat{B}_{i+1}$ such that
$$\hat{h}_{i+1}(x_0)+\delta\hat{h}_{i-1}(x_0)=0.$$
By the above discussions and \eqref{zero}, $\hat{h}_{i}(x_0)$ must
be a zero point of $S_{p+1}|_{[-2,2]}$, then there is a unique $1\le
l\le p$ such that
$$\hat{h}_{i}(x_0)=2\cos\frac{l\pi}{p+1}.$$
For any $y_j\in \hat{B}_{i}$ with
$\hat{h}_{i}(y_j)=2\cos\frac{j\pi}{p}$, $j=1,\cdots,p-1$, by
\eqref{ladder-i} and a simple computation, we get
$$|\hat{h}_{i+1}(y_j)|=|z_\pm(\hat{h}_i(y_j),\hat{h}_{i-1}(y_j),V)|\ge V-2,$$
which yields $\hat{h}_{i}(\hat{B}_{i+1})\subset
2\cos[\frac{l-1}{p}\pi,\frac{l}{p}\pi]$. Hence, for any $x\in
\hat{B}_{i+1}$, there exists a unique $c$ with $|c|<1$ such that
$\hat{h}_{i}(x)=2\cos\frac{l+c}{p+1}\pi$.

For any $x\in \hat{B}_{i+1}$, by $|S_p(2\cos\theta)|\leq |S_{p+1}(2\cos\theta)|+1$,
$$\begin{array}{rcl}
2&\geq& |\hat{h}_{i+1}(x)|\\
&\geq& |z_\pm(x) S_{p+1}(\hat{h}_{i}(x))|-
|\hat{h}_{i-1}(x)S_p(\hat{h}_{i}(x))|\\
&\geq& (V-2)|S_{p+1}(\hat{h}_{i}(x))|-
(|S_{p+1}(\hat{h}_{i}(x))|+1) |\hat{h}_{i-1}(x)|\\
&=& (V-4)|S_{p+1}(\hat{h}_{i}(x))|-2,
\end{array}$$
which follows that $|S_{p+1}(\hat{h}_{i}(x))|\leq\frac{4}{V-4}$.
Hence if $V>20$, we have $|S_{p+1}(\hat{h}_{i}(x))|\leq\frac{1}{4}$.
A direct computation gives also
$$S_{p+1}(\hat{h}_{i}(x))=S_{p+1}(2\cos\frac{(l+c)\pi}{p+1})
=\frac{\sin(l+c)\pi}{\sin\frac{l+c}{p+1}\pi} =\frac{(-1)^l\sin
c\pi}{\sin\frac{l+c}{p+1}\pi}.$$ we get finally
$|c|\le\frac{1}{10}$.
\end{proof}

\begin{defi}
Assume $V>20$. Let $(\hat{B}_i)_{i=0}^m$ be a modified ladder. Let
$(p_i)_{i=1}^{m-1}$ and $(l_i)_{i=1}^{m-1}$ be given as in
\eqref{pvalue} and Proposition A\ref{index} respectively, which will
be called the {\em type sequence} and the {\em index sequence}
respectively with respect to the modified ladder.
\end{defi}

Note that if $\alpha=[0; a_1,a_2,\cdots]$ with $(a_n)_{n\ge1}$ bounded by a constant $M$,
then the type sequence $(p_i)_{i=1}^{m-1}$ is bounded by $M+1$.

\begin{llem}\label{rit}
{\rm (\cite{R,LW})} For $V>4$, $k\ge0$,

$(1)$ A $(k,{\rm I})$-type band contains a unique band of
$\sigma_{(k+2,0)}$ which is a $(k+1,{\rm II})$-type band;

$(2)$ Let $B_{(k+1,0)}$ be a $(k,{\rm II})$-type band.

$B_{(k+1,0)}$ contains $a_{k+1}+1$ bands of $\sigma_{(k+1,1)}$ which
are of $(k+1,{\rm I})$-type, note that the fact
$$t_{(k+1,0)}({B}_{(k+1,1)}^{(i)})\subset I_{a_{k+1}+1,i}\ ,\quad
i=1,\cdots,a_{k+1}+1,$$ we can index these bands as
$\{{B}_{(k+1,1)}^{(i)}\}_{i=1}^{a_{k+1}+1}$.

$B_{(k+1,0)}$ contains $a_{k+1}$ bands of $\sigma_{(k+2,0)}$, which
are of $(k+1,{\rm III})$-type, and we can index them as
$\{B_{(k+2,0)}^{(i)}\}_{i=1}^{a_{k+1}}$ by the fact
$$t_{(k+1,0)}({B}_{(k+2,0)}^{(i)})\subset I_{a_{k+1},i}\ ,\quad
i=1,\cdots,a_{k+1}.$$

$(3)$ Let $B_{(k+1,0)}$ be a $(k,{\rm III})$-type band.

$B_{(k+1,0)}$ contains $a_{k+1}$ bands
of $\sigma_{(k+1,1)}$, which are of $(k+1,{\rm I})$-type,
and we can index them as $\{{B}_{(k+1,1)}^{(i)}\}_{i=1}^{a_{k+1}}$ by the fact
$$t_{(k+1,0)}({B}_{(k+1,1)}^{(i)})\subset I_{a_{k+1},i}\ ,\quad
i=1,\cdots,a_{k+1}.$$

$B_{(k+1,0)}$ contains $a_{k+1}-1$ bands of $\sigma_{(k+2,0)}$,
which are of $(k+1,{\rm III})$-type, and we can index them
as $\{B_{(k+2,0)}^{(i)}\}_{i=1}^{a_{k+1}-1}$ by the fact
$$t_{(k+1,0)}({B}_{(k+2,0)}^{(i)})\subset I_{a_{k+1}-1,i}\ ,\quad
i=1,\cdots,a_{k+1}-1.$$
\end{llem}

We summarize the
estimation of Chebischev polynomials on the interval $I_{p,l}$,
which has been got in the proof of Proposition 7 of \cite{LW}.

\begin{lprop}\label{keyLW}
Fix $p\ge1$, $1\le l\le p$. For $V>20$, and any $t\in I_{p,l}$,
$$\begin{array}{l}
|S_{p+1}(t)|\le\frac{1}{4},\quad |S_p(t)|\le \frac{5}{4},\\
\frac{p+1}{3}\le |S'_{p+1}(t)|\le \frac{(p+1)^3}{4},\quad |S'_p(t)|\le2|S'_{p+1}(t)|.
\end{array}
$$
\end{lprop}

With above discussions, we can simplify part of the statement of
Proposition 7,8,9 of \cite{LW} as the following, which is got by
Proposition A\ref{keyLW} and \eqref{ladder-i}.

\begin{lprop}[\cite{LW}]\label{lm-2}
Assume $V>20$. Let $(\hat{B}_i)_{i=0}^m$ be a modified ladder,
$(\hat{h}_i)_{i=0}^m$ the corresponding generating polynomials and
$(p_i)_{i=1}^{m-1}$ the corresponding type sequence. For any
$0<i<m$, we have,
$$(V-8)\frac{(p_i+1)}{3}\le\left|\frac{\hat{h}'_{i+1}(x)}{\hat{h}'_{i}(x)}\right|
\le (V+8)\frac{(p_i+1)^3}{4}.$$
\end{lprop}

Now we state the results of this paper which are motivated from
\cite{MRW}.

Fix $\alpha=[0;a_1,a_2,a_3,\cdots]$ with $(a_n)_{n\ge1}$ bounded by $M(\ge1)$.

\begin{thm}[Bounded variation]\label{bvar}
Let $V>20$, $B_n$ be a spectral generating band of order $n$ with generating polynomial $h_n$.
There exists a constant $\xi\ge1$ independent of $n$ and $B_n$ such that,
for any $x_1,x_2\in B_n$,
$$\xi^{-1}\le \left|\frac{h'_n(x_1)}{h'_n(x_2)}\right|\le \xi.$$
\end{thm}

\begin{cor}[Bounded distortion]\label{bdist}
Let $V>20$, $B_n$ be a spectral generating band of order $n$ with generating polynomial $h_n$.
Then for any $x\in B_n$,
$$\xi^{-1}\le|h'_n(x)|\cdot|B_n|\le \xi.$$
\end{cor}

\begin{thm}[Bounded covariation]\label{bco}
Suppose $V$ is sufficiently large. Given $n>k\ge1$, let
$$\begin{array}{l}
B_{n}\subset \cdots\subset B_{k+1}\subset B_{k}\\
\tilde{B}_{n}\subset\cdots\subset \tilde{B}_{k+1}\subset \tilde{B}_{k},
\end{array}
$$
be two sequence of spectral generating bands. For any $k+1\le i\le
n$, $B_i$ and $\tilde{B}_i$ are of same type and have the same
index. $B_k$ and $\tilde{B}_k$ are of same type.
So, there exists $\eta>1$ such that
$$\eta^{-1} \frac{|\tilde{B}_n|}{|\tilde{B}_{k}|}\le
\frac{|B_n|}{|B_{k}|}\le \eta \frac{|\tilde{B}_n|}{|\tilde{B}_{k}|}.$$
\end{thm}

\begin{thm}[Existence of Gibbs-like measures]\label{gibbs}
Suppose $V$ is sufficiently large. Given $\beta>0$,
there exist $\zeta>0$ and a probability measure $\mu_\beta$
supported by $\sigma(H_{\alpha,V})$ such that for any $k\ge1$ and $\tilde{B}\in \mathscr{G}_k$,
$$\zeta^{-1}\frac{|\tilde{B}|^\beta}{\sum\limits_{B\in\mathscr{G}_k}|B|^\beta}
\le \mu_\beta(\tilde{B})\le \zeta
\frac{|\tilde{B}|^\beta}{\sum\limits_{B\in\mathscr{G}_k}|B|^\beta}.
$$
\end{thm}

Let $s_n$ be the $n$-th pre-dimension of $\sigma(H_{\alpha,V})$, i.e.,
$$\sum_{B\in\mathscr{G}_n}|B|^{s_n}=1,$$
and
$$s_*=\liminf_{n\rightarrow\infty}s_n,\quad s^*=\limsup_{n\rightarrow\infty}s_n.$$

\begin{thm}\label{haus}
For sufficiently large $V$, $\dim_H\, \sigma(H_{\alpha,V})=s_*$.
\end{thm}

\noindent {\bf Remark:} Liu, Peyri\`ere and Wen proved in
\cite{LPW07} that $\overline{\dim}_B\, \sigma(H_{\alpha,V})=s^*$,
but there was a small error there (Page 670 in \cite{LPW07}), we can
correct it as follows:

letting $B_n\subset B_{n-1}$ be two spectral generating bands of
order $n$ and $n-1$, taking notation of \cite{LPW07}, denoting $B_n$
as $J$, $B_{n-1}$ as $J^{-1}$, the inequality
\begin{equation*}\label{error}
q_{ij}(n)\le \frac{|J|}{|J^{-1}|}\le p_{ij}(n)
\end{equation*}
should be replaced by the inequality(see \cite{LW}, in proof of Proposition 5, $(4.36)$-$(4.40)$)
$$q_{ij}(n)\le \frac{|h_{n-1}'(x)|}{|h_n'(x)|}\le p_{ij}(n),\quad \forall x\in J,$$
where $h_n$ and $h_{n-1}$ are the corresponding generating
polynomial of $B_n$ and $B_{n-1}$.

From above inequality, applying Corollary $\ref{bdist}$ we get
$$\frac{q_{ij}(n)}{\xi^2}\le \frac{|J|}{|J^{-1}|}\le \min\{\xi^2 p_{ij}(n),1\},$$
which yields the lower bound of contractive ratio is strictly larger
than $0$. Then as in the proof of \cite{LPW07}, we still have

\begin{lthm}For $V>20$, $\overline{\dim}_B\, \sigma(H_{\alpha,V})=s^*$.
\end{lthm}

\mysection{Bounded variation, Bounded covariation, Gibbs measure}
%
%As we have showed, the key step for estimating the dimensions of the
%spectrum is to estimate the length of the spectral generating bands,
%which is reduced to estimate the derivative of the generating
%polynomial on the different points of correspondent spectral
%generating bands. The sketch of the method in \cite{LW} can be
%summarize in the following way.
%
%For any $n>1$ and $B_n\in\mathscr{G}_n$, there exists a unique
%sequence $(B_i)_{-1\le i<n}$ such that $B_i\in\mathscr{G}_i$ and
%$$B_n\subset B_{n-1}\subset\cdots\subset B_0\subset B_{-1}.$$
%Let $(h_i)_{-1\le i\le n}$ be the corresponding generating
%polynomial, for any $x\in B_n$ and $1\le i<n$, we will estimate
%$\frac{h_i'(x)}{h_{i+1}'(x)}$. If $B_i$ is either $(i,II)$-type or
%$(i,III)$-type band, we can estimate it directly by \eqref{cII} or
%\eqref{cIII}. Now we consider $B_i$ is of $(i,I)$-type, in this
%case, the estimatings depend on $a_i$. If $a_i=1$, then
%$B_i=B_{i+1}$ and $h_i=h_{i+1}$ by the facts $t_{(i+2,0)}=t_{(i,1)}$
%and $t_{(i+1,-1)}=t_{(i,0)}$; if $a_i>1$, then $B_{i+1}$ is a band
%in $\sigma_{(i,a_i)}$, there exists a unique sequence of bands such
%that
%$$B_{i+1}=B_{(i,a_i)}\subset B_{(i,a_i-1)}\subset\cdots\subset
%B_{(i,1)}=B_i.$$ Using \eqref{cI},  for $1\le p<a_i$, we can
%estimate $\frac{t'_{(i,p)}}{t'_{(i,p+1)}}$, then by the relation
%above, we get successively
%$$\frac{h_i'}{h_{i+1}'}=\frac{t'_{(i,1)}}{t'_{(i,2)}}
%\frac{t'_{(i,2)}}{t'_{(i,3)}}\cdots
%\frac{t'_{(i,a_i-1)}}{t'_{(i,a_i)}}.$$

From Proposition A\ref{lm-2}, we get immediately the following corollary.

\begin{cor}\label{expand}
Assume $V>20$. Let $(\hat{B}_i)_{i=0}^m$ be a modified ladder and
$(\hat{h}_i)_{i=0}^m$ the correspondent generating polynomials. Then
for any $x,y\in \hat{B}_m$,
$$|\hat{h}_i(x)-\hat{h}_i(y)|\le 3^{-(m-i)}
|\hat{h}_m(x)-\hat{h}_m(y)|\le 4\cdot 3^{-(m-i)}.$$
\end{cor}
\begin{proof}
For any $0\le i\le m$, since $\hat{h}_i$ is monotone on $\hat{B}_i$
$$\begin{array}{rcl}
|\hat{h}_i(x)-\hat{h}_i(y)|&=&\left|\int_x^y \hat{h}'_i(t) dt\right|\\
&=&\left|\int_x^y \frac{\hat{h}'_i(t)}{\hat{h}'_{i+1}(t)} \hat{h}'_{i+1}(t) dt\right|\\
&\le& 3^{-1} \left|\int_x^y \hat{h}'_{i+1}(t) dt\right|\\
&=& 3^{-1} |\hat{h}_{i+1}(x)-\hat{h}_{i+1}(y)|,\\
\end{array}$$
where the inequality is due to Proposition A\ref{lm-2}.
\end{proof}

\begin{prop}\label{lip} Assume $V>20$. Let $(\hat{B}_i)_{i=0}^m$ and
$(\hat{\tilde{B}}_i)_{i=0}^m$ be two modified ladders. Suppose that
they have the same sequence of generating polynomials
$(\hat{h}_i)_{i=0}^m$ and the same type sequence
$(p_i)_{i=1}^{m-1}$.

Suppose $x_1\in \hat{B}_{m}$, $x_2\in \hat{\tilde{B}}_{m}$, $0<i<m$.

In the case of $\hat{B}_i$ is not a band of order $0$, then
\begin{equation}\label{prop2}
\begin{array}{rcl}
\left|\dfrac{\hat{h}'_{i+1}(x_1)}{h'_{i}(x_1)}-
\dfrac{\hat{h}'_{i+1}(x_2)}{\hat{h}'_{i}(x_2)}\right|
&\le& C(|\hat{h}_i(x_1)-\hat{h}_i(x_2)|+
|\hat{h}_{i-1}(x_1)-\hat{h}_{i-1}(x_2)|)+\\
&&\ \dfrac{1}{3}
\left|\dfrac{\hat{h}'_i(x_1)}{\hat{h}'_{i-1}(x_1)}-
\dfrac{\hat{h}'_i(x_2)}{\hat{h}'_{i-1}(x_2)}\right|,
\end{array}
\end{equation}
where $C$ is a constant depending on $V$ and $p_i$.

In the case of $\hat{B}_i$ is a band of order $0$, then
\begin{equation}\label{prop2-0}
\begin{array}{rcl}
\left|\dfrac{\hat{h}'_{i+1}(x_1)}{h'_{i}(x_1)}-
\dfrac{\hat{h}'_{i+1}(x_2)}{\hat{h}'_{i}(x_2)}\right| &\le&
C|\hat{h}_i(x_1)-\hat{h}_i(x_2)|.
\end{array}
\end{equation}
\end{prop}
\begin{proof}
%Fix any $x_1\in \hat{B}_{m}$, $x_2\in \hat{\tilde{B}}_{m}$, $1<i<m$,
%we will discuss all possible types of $\hat{B}_i$, that is,
%$\hat{B}_i$ may be one of the following five cases: type $I$, $II$,
%$III$ of order $k$, an added rung of order $k-1$ and an added rung
%of order $k$.
%
%Case 1)\ \ If $\hat{B}_i$ is $(k,II)$-type, then
%$\hat{h}_i=t_{(k+1,0)}$, $\hat{h}_{i+1}=t_{(k,p)}$ for some $p\ge1$,
%$\hat{h}_{i-1}=t_{(k,-1)}=t_{(k-1,a_k-1)}$ (note that
%$\hat{B}_{i-1}$ is an added rung if $a_{k}>2$, this is also an
%advantage to apply ladder ). We have
%\begin{equation}\label{kii}
%\begin{array}{rcl}
%\hat{h}_{i+1}=t_{(k,p)}&=&t_{(k,0)}S_{p+1}(t_{(k+1,0)})-t_{(k,-1)}S_p(t_{(k+1,0)})\\
%&=&t_{(k,0)}S_{p+1}(\hat{h}_{i})-\hat{h}_{i-1}S_p(\hat{h}_i)\\
%t_{(k,0)}&=&z_{\pm}(\hat{h}_{i},\hat{h}_{i-1},V).
%\end{array}\end{equation}
Take any $0<i<m$, for convenience, we denote
$z_{\pm}(\hat{h}_i(x),\hat{h}_{i-1}(x),V)$ as $z_\pm(x)$.

Suppose first $\hat{B}_i$ is not a band of order $0$.

%By \eqref{ladder-i}, we have
%$$\hat{h}_{i+1}(x)=z_{\pm}(x)
%S_{p+1}(\hat{h}_{i}(x))-\hat{h}_{i-1}(x)S_p(\hat{h}_i(x)),$$
%where $p$ is given according \eqref{pvalue}.
%
%We claim that
%there exist $l(1\le l\le p)$ such that
%for any $x\in \hat{B}_{i+1}\cup
%\hat{\tilde{B}}_{i+1}$,
%$$\hat{h}_i(x)\in I_{p,l}.$$
%If $\hat{B}_i=B_k$ is a band of $(k,II)$ or $(k,III)$-type,
%then this is get by hypothesis of $B_{k+1}$ and $\tilde{B}_{k+1}$ are of same type and index.
%If $\hat{B}_i$ is a band of type $I$ or an added rung, then by $p=1$
%and Proposition \ref{index}, the claim holds.

By taking derivative on both side of \eqref{ladder-i}, we get
\begin{equation}\label{exprderi}
\begin{array}{r}
\dfrac{\hat{h}'_{i+1}(x)}{\hat{h}'_i(x)}
=S'_{p_i+1}(\hat{h}_i(x))z_\pm(x)-S'_{p_i}(\hat{h}_i(x))\hat{h}_{i-1}(x)+\quad\\
S_{p_i+1}(\hat{h}_i(x))\dfrac{z'_\pm(x)}{\hat{h}'_i(x)}-
S_{p_i}(\hat{h}_i(x))\dfrac{\hat{h}'_{i-1}(x)}{\hat{h}'_i(x)}.
\end{array}
\end{equation}

Observing that $S_{p+1}(x)$ is a polynomial of degree $p$,
$S''_{p+1}|_{[-2,2]}$ is also bounded by some constant depends on
$p$. So, there exists a constant $c_1>0$ depending only on $p_i$
such that
\begin{equation}\label{cha1}
\begin{array}{l}
\left|S_{p_i+1}(\hat{h}_i(x_1))-S_{p_i+1}(\hat{h}_i(x_2))\right|
\le c_1 |\hat{h}_i(x_1)-\hat{h}_i(x_2)|\\[8pt]
\left|S'_{p_i+1}(\hat{h}_i(x_1))-S'_{p_i+1}(\hat{h}_i(x_2))\right|
\le c_1 |\hat{h}_i(x_1)-\hat{h}_i(x_2)|\\[8pt]
\left|S_{p_i}(\hat{h}_i(x_1))-S_{p_i}(\hat{h}_i(x_2))\right|
\le c_1 |\hat{h}_i(x_1)-\hat{h}_i(x_2)|\\[8pt]
\left|S'_{p_i}(\hat{h}_i(x_1))-S'_{p_i}(\hat{h}_i(x_2))\right|
\le c_1 |\hat{h}_i(x_1)-\hat{h}_i(x_2)|\\[8pt]
\left|z_\pm(x_1)-z_\pm(x_2)\right|
\le |\hat{h}_i(x_1)-\hat{h}_i(x_2)|+|\hat{h}_{i-1}(x_1)-\hat{h}_{i-1}(x_2)|\\
\end{array}
\end{equation}
where the last inequality is due to the fact \eqref{dec}, and
\begin{equation}\label{d-cha}
\frac{z'_\pm(x)}{\hat{h}'_i(x)}=
z_1(\hat{h}_i(x),\hat{h}_{i-1}(x),V)+z_2(\hat{h}_i(x),\hat{h}_{i-1}(x),V)
\frac{\hat{h}'_{i-1}(x)}{\hat{h}'_i(x)}.
\end{equation}

By \eqref{dec}, we have
%\begin{equation}\label{cha2}
%\begin{array}{cl}
%&|z_1(\hat{h}_i(x_1),\hat{h}_{i-1}(x_1),V)-z_1(\hat{h}_i(x_2),\hat{h}_{i-1}(x_2),V)|\\
%\le&
%|\hat{h}_i(x_1)-\hat{h}_i(x_2)|+|\hat{h}_{i-1}(x_1)-\hat{h}_{i-1}(x_2)|\\
%&|z_2(\hat{h}_i(x_1),\hat{h}_{i-1}(x_1),V)-z_2(\hat{h}_i(x_2),\hat{h}_{i-1}(x_2),V)|\\
%\le&
%|\hat{h}_i(x_1)-\hat{h}_i(x_2)|+|\hat{h}_{i-1}(x_1)-\hat{h}_{i-1}(x_2)|.\\
%\end{array}\end{equation}
\begin{equation}\label{cha2}
\begin{array}{l}
|z_1(\hat{h}_i(x_1),\hat{h}_{i-1}(x_1),V)-z_1(\hat{h}_i(x_2),\hat{h}_{i-1}(x_2),V)|\\
\quad\le |\hat{h}_i(x_1)-\hat{h}_i(x_2)|+|\hat{h}_{i-1}(x_1)-\hat{h}_{i-1}(x_2)|\\[5pt]
|z_2(\hat{h}_i(x_1),\hat{h}_{i-1}(x_1),V)-z_2(\hat{h}_i(x_2),\hat{h}_{i-1}(x_2),V)|\\
\quad\le |\hat{h}_i(x_1)-\hat{h}_i(x_2)|+|\hat{h}_{i-1}(x_1)-\hat{h}_{i-1}(x_2)|.\\
\end{array}\end{equation}
By a direct computation,
\begin{equation}\label{cha3}
\begin{array}{rcl}
\left|\dfrac{\hat{h}'_{i-1}(x_1)}{\hat{h}'_i(x_1)}-
\dfrac{\hat{h}'_{i-1}(x_2)}{\hat{h}'_i(x_2)}\right|
&=&\left|\dfrac{\hat{h}'_{i-1}(x_1)}{\hat{h}'_i(x_1)}
\dfrac{\hat{h}'_{i-1}(x_2)}{\hat{h}'_i(x_2)}\right|
\left|\dfrac{\hat{h}'_{i}(x_1)}{\hat{h}'_{i-1}(x_1)}-
\dfrac{\hat{h}'_{i}(x_2)}{\hat{h}'_{i-1}(x_2)}\right|\\
&\le&\dfrac{1}{9}\left|\dfrac{\hat{h}'_{i}(x_1)}{\hat{h}'_{i-1}(x_1)}-
\dfrac{\hat{h}'_{i}(x_2)}{\hat{h}'_{i-1}(x_2)}\right|.
\end{array}
\end{equation}

The inequalities \eqref{exprderi}-\eqref{cha3} imply that the
inequality \eqref{prop2} holds.

%Since there is no bands of $(0,II)$-type, we need only consider
%band.
\medskip

Suppose $\hat{B}_i$ is a band of order $0$. Note that
$\hat{h}_{i-1}=t_{(0,0)}\equiv2$ is a constant, then an analogous
argument to \eqref{exprderi}-\eqref{cha2} implies that the
inequality \eqref{prop2-0} holds.
\end{proof}

\begin{proof}[Proof of Theorem \ref{bvar}]
It is a corollary of Corollary \ref{expand} and Proposition
\ref{lip}. In fact, let $$B_n\subset B_{n-1}\subset\cdots\subset
B_0\subset B_{-1}$$ be a sequence of spectral generating bands(the
orders are from $n$ to $-1$), which form an initial ladder. Let
$(\hat{B}_i)_{i=-1}^m$ be the corresponding modified ladder,
$(\hat{h}_i)_{i=-1}^m$ the corresponding generating polynomials. By
Corollary \ref{expand}, for any $0\le i<n$,
$$|\hat{h}_i(x_1)-\hat{h}_i(x_2)|\le 4\cdot 3^{-(m-i)}.$$

Note that $\hat{B}_0=B_0$, $\hat{B}_{-1}={B}_{-1}$, we have
$\hat{h}'_0\equiv1$, thus
\begin{equation}\label{thm1-1}
|\log|\hat{h}'_m(x_1)|-\log|\hat{h}'_m(x_2)||
\leq\sum\limits_{i=0}^{m-1}\left|
\log\left|\dfrac{\hat{h}'_{i+1}(x_1)}{\hat{h}'_{i}(x_1)}\right|-
\log\left|\dfrac{\hat{h}'_{i+1}(x_2)}{\hat{h}'_{i}(x_2)}\right|\right|
\end{equation}
and
\begin{equation}\label{thm1-2}
\left|\log\left|\dfrac{\hat{h}'_{i+1}(x_1)}{\hat{h}'_{i}(x_1)}\right|-
\log\left|\dfrac{\hat{h}'_{i+1}(x_2)}{\hat{h}'_{i}(x_2)}\right|\right|
\le\left|\dfrac{\hat{h}'_{i}(x_2)}{\hat{h}'_{i+1}(x_2)}\right|
\left|\dfrac{\hat{h}'_{i+1}(x_1)}{\hat{h}'_{i}(x_1)}-
\dfrac{\hat{h}'_{i+1}(x_2)}{\hat{h}'_{i}(x_2)}\right|.
\end{equation}

$\hat{B}_0$ is a order $0$ band, so by \eqref{prop2-0}
$$
\left|\dfrac{\hat{h}'_{1}(x_1)}{\hat{h}'_{0}(x_1)}-
\dfrac{\hat{h}'_{1}(x_2)}{\hat{h}'_{0}(x_2)}\right|\le C |\hat{h}_0(x_1)-\hat{h}_0(x_2)|
\le 4C \cdot 3^{-m}.
$$
Combining with \eqref{prop2} and induction, we have for $0\le i<m$,
$$
\left|\dfrac{\hat{h}'_{i+1}(x_1)}{\hat{h}'_{i}(x_1)}-
\dfrac{\hat{h}'_{i+1}(x_2)}{\hat{h}'_{i}(x_2)}\right|
\le 8C \cdot 3^{-(m-i)},
$$
together with \eqref{thm1-1} and \eqref{thm1-2}, we finish the proof
of the theorem.
\end{proof}

\begin{prop}\label{diff}
Suppose that $(\hat{B}_i)_{i=0}^m$ and $(\hat{\tilde{B}}_i)_{i=0}^m$
are two modified ladders having the same sequence of generating
polynomials $(\hat{h}_i)_{i=0}^m$, the same type sequence
$(p_i)_{i=1}^{m-1}$(bounded by $M+1$), and the same index sequence
$(l_i)_{i=1}^{m-1}$. Then for sufficiently large $V$, we have
\begin{itemize}
\item[{\rm(i)}]\ There exists $c<\frac{1}{4}$ such that
for any $0<i<m$ and any $x_1\in \hat{B}_{i+1}$, $x_2\in \hat{\tilde{B}}_{i+1}$,
\begin{equation}\label{ediff}
|\hat{h}_{i}(x_1)-\hat{h}_i(x_2)|\le c
(|\hat{h}_{i+1}(x_1)-\hat{h}_{i+1}(x_2)|+|\hat{h}_{i-1}(x_1)-\hat{h}_{i-1}(x_2)|).
\end{equation}
\item[{\rm (ii)}]\ Letting $\lambda=\frac{1+\sqrt{1-4c^2}}{2c}(>1)$,
for any $x_1\in \hat{B}_{m}$, there exists $x_2\in \hat{\tilde{B}}_{m}$ such that,
for any $0\le i\le m$,
\begin{equation}\label{distdist}
|\hat{h}_{i}(x_1)-\hat{h}_i(x_2)|\le \frac{4\lambda^2}{\lambda^2-1} \lambda^{-i}.
\end{equation}
\item[{\rm (iii)}]\ There exists $\eta>1$ such that
$$\eta^{-1} \frac{|\hat{\tilde{B}}_m|}{|\hat{\tilde{B}}_{0}|}\le
\frac{|\hat{B}_m|}{|\hat{B}_{0}|}\le \eta \frac{|\hat{\tilde{B}}_m|}{|\hat{\tilde{B}}_{0}|}.$$
\end{itemize}
\end{prop}

\begin{proof}
(i) Take any $0<i<m$ and any $x_1\in \hat{B}_{i+1}$, $x_2\in
\hat{\tilde{B}}_{i+1}$. For convenience, we denote
$z_{\pm}(\hat{h}_i(x),\hat{h}_{i-1}(x),V)$ by $z_\pm(x)$.

By the definitions of $p_i$ and $l_i$, $\hat{h}_i(x_1)$,
$\hat{h}_i(x_2)\in I_{p_i,l_i}$, then by Proposition A\ref{keyLW},
\begin{equation}\label{indexkey-2}
\left|S_{p_i+1}(\hat{h}_i(x_1))-S_{p_i+1}(\hat{h}_i(x_2))\right|
\ge \frac{p_i+1}{3} |\hat{h}_i(x_1)-\hat{h}_i(x_2)|.
\end{equation}

By Proposition A\ref{keyLW} again,
$$\left|S_{p}(\hat{h}_i(x_1))-S_{p}(\hat{h}_i(x_2))\right|
\le \frac{(p+1)^3}{4} |\hat{h}_i(x_1)-\hat{h}_i(x_2)|.$$
By \eqref{dec} as in \eqref{cha1},
$$
\left|z_\pm(x_1))-z_\pm(x_2))\right| \le
|\hat{h}_i(x_1)-\hat{h}_i(x_2)|+|\hat{h}_{i-1}(x_1)-\hat{h}_{i-1}(x_2)|.
$$
So by the above three inequalities and \eqref{ladder-i}, we get
$$\begin{array}{cl}
|\hat{h}_{i+1}(x_1)-\hat{h}_{i+1}(x_2)|
\ge& \frac{(V-2)(p+1)}{3}|\hat{h}_i(x_1)-\hat{h}_i(x_2)|
-\frac{1}{4}|\hat{h}_i(x_1)-\hat{h}_i(x_2)|\\[4pt]
& -\frac{6}{4}|\hat{h}_{i-1}(x_1)-\hat{h}_{i-1}(x_2)|-
{(p+1)^3}|\hat{h}_i(x_1)-\hat{h}_i(x_2)|,\\
\end{array}$$
which concludes the inequality \eqref{ediff}.

\medskip

(ii) Since
$$%\begin{array}{l}
\hat{h}_m(\hat{B}_m)=[-2,2],\ \hat{h}_m(\hat{\tilde{B}}_m)=[-2,2],
%\end{array}
$$
for any $x_1\in \hat{B}_m$, there exists $x_2\in\hat{\tilde{B}}_m$
such that
$$\hat{h}_m(x_1)=\hat{h}_m(x_2).$$

Since for any $0\le i<m$,
$$\hat{h}_i(\hat{B}_m)\subset[-2,2],\
\hat{h}_i(\hat{\tilde{B}}_m)\subset[-2,2],$$ we get for any $0\le
i<m$,
$$|\hat{h}_i(x_1)-\hat{h}_i(x_2)|\le4.$$

Let $$f_i=|\hat{h}_{i}(x_1)-\hat{h}_i(x_2)|,\quad 0\le i\le m,$$
then $f_m=0$, $f_0\le4$. By \eqref{ediff},
$$0\le(\lambda f_{m-1}-f_m)\le \lambda^{-1}(\lambda f_{m-2}-f_{m-1})
\le\cdots\le
\lambda^{-m+1}(\lambda f_0-f_{1})\le 4\lambda^{-m+2},$$
which implies that for any $0\le i\le m$
$$|\hat{h}_{i}(x_1)-\hat{h}_i(x_2)|=f_i\le\frac{4\lambda^2}{\lambda^2-1}\lambda^{-i}.
$$

(iii) By \eqref{distdist}, an argument similar to Theorem \ref{bvar} implies
there exist $\xi_1>1$ such that
\begin{equation}\label{ratio}
\xi_1^{-1}\le\left|\frac{\hat{h}'_m(x_1)/\hat{h}'_0(x_1)}
{\hat{h}'_m(x_2)/\hat{h}'_0(x_2)}\right|\le \xi_1.
\end{equation}

By the definition of the generating polynomial, there exist
$\tilde{x}\in \hat{B}_m$, $\tilde{y}\in\hat{B}_0$ such that
$$|\hat{B}_m|\,|\hat{h}_m'(\tilde{x})|=4,\ |\hat{B}_0|\,|\hat{h}_0'(\tilde{y})|=4.$$

Associating with Theorem \ref{bvar}, we have
$$\frac{|\hat{B}_m|}{|\hat{B}_0|}=
\frac{|\hat{B}_m|\,|\hat{h}_m'(\tilde{x})|}{|\hat{B}_0|\,|\hat{h}_0'(\tilde{y})|}
\left|\frac{\hat{h}_m(x_1)}{\hat{h}_m(\tilde{x})}\right|\,
\left|\frac{\hat{h}_0(\tilde{y})}{\hat{h}_0(x_1)}\right|\,
\left|\frac{\hat{h}_0(x_1)}{\hat{h}_m(x_1)}\right|
\le \xi^2\left|\frac{\hat{h}_0(x_1)}{\hat{h}_m(x_1)}\right|.
$$

By the same discussion, we have
$$\frac{|\hat{\tilde{B}}_m|}{|\hat{\tilde{B}}_0|}
\ge \xi^{-2}\left|\frac{\hat{h}_0(x_2)}{\hat{h}_m(x_2)}\right|.
$$

Then by \eqref{ratio}, we have
$$\frac{|\hat{B}_m|}{|\hat{B}_0|}\le \xi^4\xi_1\frac{|\hat{\tilde{B}}_m|}{|\hat{\tilde{B}}_0|}.$$

The opposite direction of the inequality can be got by the same way.
\end{proof}

\begin{proof}[Proof of Theorem \ref{bco}]
Let $(\hat{B}_i)_{i=0}^m$ be the modified ladder of initial ladder
$({B}_i)_{i=k}^n$ and $(\hat{\tilde{B}}_i)_{i=0}^m$ be the modified
ladder of the initial ladder $(\tilde{B}_i)_{i=k}^n$.

Since for $k\le i\le m$, $B_i$ and $\tilde{B}_i$ are of the same
type, $(\hat{B}_i)_{i=0}^m$ and $(\hat{\tilde{B}}_i)_{i=0}^m$ share
the same sequence of generating polynomials $(\hat{h}_i)_{i=0}^m$
and the same type sequence $(p_i)_{i=1}^{m-1}$.

Since for $k< i\le m$, $B_i$ and $\tilde{B}_i$ are of the same
index, $(\hat{B}_i)_{i=0}^m$ and $(\hat{\tilde{B}}_i)_{i=0}^m$ share
the same index sequence $(l_i)_{i=1}^{m-1}$.

Then Proposition \ref{diff} concludes the result of the theorem.
\end{proof}

Given a spectral generating band $B_n$ of order $n$, there exists a
unique sequence of spectral generating bands $(B_i)_{i=0}^{n-1}$ so
that
$$B_n\subset B_{n-1}\subset \cdots B_1\subset B_0,$$
we are going to define the {\em characteristic index $i_0 i_1 \cdots
i_n$} of $B_n$ as follows, fix $0\le k\le n-1$,
\begin{itemize}
\item[Case 1:]\ $B_k$ is a $(k,II)$-type band.\\
If $B_{k+1}$ is $(k+1,I)$ type band with index $j$, then define $i_{k+1}:=(I,j)$;\\
if $B_{k+1}$ is $(k+1,III)$ type band with index $j$, then define
$i_{k+1}:=(III,j)$.

\item[Case 2:]\ $B_k$ is a $(k,III)$-type band.\\
If $B_{k+1}$ is $(k+1,I)$ type band with index $j$, then $i_{k+1}:=(I,j)$;\\
if $B_{k+1}$ is $(k+1,III)$ type band with index $j$, then
$i_{k+1}:=(III,j)$.

\item[Case 3:]\ $B_k$ is a $(k,I)$-type band, then $i_{k+1}:=(II)$.

\item[Case 4:]\ If $B_0$ is of $(0,I)$-type, then $i_0:=(I)$;
if $B_0$ is of $(0,III)$-type, then $i_0=(III)$.
\end{itemize}

We call $i_0 i_1\cdots i_m$ an {\em admissible index} (of length
$m$) if it is a characteristic index of a band $B_m$ of order $m$.
Denote by $\Omega_m$ the set of all admissible index of length $m$.
For any admissible index $\omega\in\Omega_m$, there is only one
associated spectral generating band, which we denoted as $B_\omega$.
For any $i_0\cdots i_m\in \Omega_m$ and any $0\le j\le m$, we call
the symbol $i_j$ is of type $I$ ($II$ or $III$), if the
corresponding band $B_{i_0\cdots i_j}\in \Omega_{j}$ is of type $I$
($II$ or $III$ respectively). Now we give some more notations:
\begin{itemize}
\item $\Omega_{k,m}$: all segments $i_{k}\cdots i_{m}$
of any admissible index $i_0 i_1\cdots i_m\in\Omega_m$, $m\ge k>0$.
\item $\Omega_{k+1,m}^{i_k}$:
all segments $i_{k+1}\cdots i_{m}$ of $i_0 i_1\cdots i_k
i_{k+1}\cdots i_m\in\Omega_m$, $i_{k}\in\Omega_{k,k}$.

Since it depends only on type of $i_k$, for the convenience, we
denoted it sometimes by $\Omega_{k+1,m}^{I}$,$\Omega_{k+1,m}^{II}$,
or $\Omega_{k+1,m}^{III}$.
\item $\Omega_{m}^{(k_1,j_1)(k_2,j_2)\cdots(k_l,j_l)}$:
all $i_0\cdots i_m\in\Omega_m$ satisfying $i_{k_s}=j_s$ with
$j_s\in\Omega_{k_s,k_s}$ for $1\le s\le l$.
\end{itemize}

For any $0<\beta<1$ and $m>0$, we define a probability
$\mu_{\beta,m}$ on $\mathbb R$  such that for any
$\omega_0\in\Omega_m$,
$$\mu_{\beta,m}(B_{\omega_0})=
\frac{|B_{\omega_0}|^\beta}{\sum\limits_{\omega\in\Omega_m}|B_{\omega}|^\beta},$$
where $\mu_{\beta,m}$ is uniformly distributed on each band
$B_\omega$. For the convenience, for any $m>0$, denote
$$b_m:=\sum_{\omega\in\Omega_m}|B_{\omega}|^\beta.$$

For any $k\ge1$, any $\omega=i_0 i_1 \cdots i_k\in \Omega_k$ and any
$m>k$, we have
$$\mu_{\beta,m}(B_{\omega})=\sum_{\sigma\in\Omega_{k+1,m}^{i_k}}\mu_{\beta,m}(B_{\omega*\sigma}),$$
where $\omega*\sigma$ is the concatenation of $\omega$ and $\sigma$.

In the following, we suppose that $V$ is large enough so that
Bounded covariation holds.

\begin{prop}\label{progib}
Let $\mu_{\beta,m}$ be defined as above. Then there exists $c\ge1$
such that
\begin{itemize}
\item[(i)]\ for any $k>0$ and $\omega\in\Omega_k$,
\begin{equation}\label{progib1}
c^{-1}\mu_{\beta,k+3}(B_{\omega})\le\mu_{\beta,k}(B_{\omega}) \le
c\mu_{\beta,k+3}(B_{\omega});
\end{equation}
\item[(ii)]\
for any $k>0$, $m>k+3$, $\omega=i_0\cdots i_k\in\Omega_k$, $\sigma\in\Omega^{i_k}_{k+1,k+3}$,
\begin{equation}\label{progib2}
\mu_{\beta,m}(B_{\omega})\le c \mu_{\beta,m}(B_{\omega*\sigma}).
\end{equation}
\end{itemize}
\end{prop}
\begin{proof}
(i) Take any $\omega_0=i_0\cdots i_k\in\Omega_k$.
For any $\sigma\in\Omega_{k+1,k+3}^{i_k}$, $x\in B_{\omega_0*\sigma}$,
by Corollary $\ref{bdist}$ and Proposition A\ref{lm-2},
$$
1\le\frac{|B_{\omega_0}|}{|B_{\omega_0*\sigma}|}\le \xi^2
\frac{|h_{k+3}'(x)|}{|h_k'(x)|}\le \xi^2((M+2)^3(V+8))^{2M+1},
$$
where $h_k(x)$ is the generating polynomial of $B_{\omega_0}$ and
$h_{k+3}(x)$ is the generating polynomial of $B_{\omega_0*\sigma}$,
and the length of modified ladder from $B_{\omega_0}$ to $B_{\omega_0*\sigma}$ is
at most $2M+1$ and at least $1$. So for
$c_1=\xi^2((M+2)^3(V+5))^{2M+1}$,
\begin{equation}\label{gibcontract}
1\le\frac{|B_{\omega_0}|}{|B_{\omega_0*\sigma}|}\le c_1.
\end{equation}

Since $B_{\omega_0}$ contains at most $(2M+1)^3$ bands of order
$k+3$,
$$(2M+1)^{-3} b_{k+3}\le b_k\le c_1 b_{k+3}.$$
Hence, for any $\omega_0\in\Omega_k$,
$$
c_1^{-1}\mu_{\beta,k+3}(B_{\omega_0}) \le\mu_{\beta,k}(B_{\omega_0})
\le c_1 (2M+1)^3\mu_{\beta,k+3}(B_{\omega_0}),
$$
which yields the inequality \eqref{progib1}.

\medskip

(ii) For any $\omega_0=i_0\cdots i_k\in \Omega_k$ and any $m>k+3$,
\begin{equation}\label{progib3}
\mu_{\beta,m}(B_{\omega_0})=\sum_{\sigma\in\Omega_{k+1,k+3}^{i_k}}
\mu_{\beta,m}(B_{\omega_0*\sigma}).
\end{equation}

We will show there exists $c_2>1$ such that, for any $m>k+6$,
$\sigma_1,\sigma_2\in \Omega_{k+1,k+3}^{i_k}$,
\begin{equation}\label{progib4}
\mu_{\beta,m}(B_{\omega_0*\sigma_1})\le c_2  \mu_{\beta,m}(B_{\omega_0*\sigma_2}).
\end{equation}
together with \eqref{progib1}, \eqref{progib3} and
$\sharp\Omega_{k+1,k+3}^{i_k}\le (2M+1)^3$, we will get the
inequality \eqref{progib2}.

Fix $\sigma_1,\sigma_2\in \Omega_{k+1,k+3}^{i_k}$. Let $i$ be the
last symbol of $\sigma_1$, $j$ be the last symbol of $\sigma_2$.
Divide $\Omega_{k+4,k+6}^{i}$ into three sets $D_1,D_2,D_3$
according to the last symbol being type $I$, type $II$, or type
$III$. Divide also the set $\Omega_{k+4,k+6}^{j}$ into three sets
$\tilde{D}_1,\tilde{D}_2,\tilde{D}_3$ according to the last symbol
being type $I$, type $II$, or type $III$. So
$$\begin{array}{l}
\mu_{\beta,m}(B_{\omega_0*\sigma_1})=
\sum\limits_{\tau\in D_1}\mu_{\beta,m}(B_{\omega_0*\sigma_1*\tau})+
\sum\limits_{\tau\in D_2}\mu_{\beta,m}(B_{\omega_0*\sigma_1*\tau})+
\sum\limits_{\tau\in D_3}\mu_{\beta,m}(B_{\omega_0*\sigma_1*\tau}),\\[10pt]
\mu_{\beta,m}(B_{\omega_0*\sigma_2})=
\sum\limits_{\tau\in \tilde{D}_1}\mu_{\beta,m}(B_{\omega_0*\sigma_2*\tau})+
\sum\limits_{\tau\in \tilde{D}_2}\mu_{\beta,m}(B_{\omega_0*\sigma_2*\tau})+
\sum\limits_{\tau\in \tilde{D}_3}\mu_{\beta,m}(B_{\omega_0*\sigma_2*\tau}),
\end{array}
$$
Fix any $\tau_1\in D_1$, $\tau_2\in \tilde{D}_1$.
$$\begin{array}{l}
\mu_{\beta,m}(B_{\omega_0*\sigma_1*\tau_1})=\sum_{\tau\in\Omega_{k+7,m}^{I}}
\mu_{\beta,m}(B_{\omega_0*\sigma_1*\tau_1*\tau}),\\[4pt]
\mu_{\beta,m}(B_{\omega_0*\sigma_2*\tau_2})=\sum_{\tau\in\Omega_{k+7,m}^{I}}
\mu_{\beta,m}(B_{\omega_0*\sigma_2*\tau_2*\tau}).
\end{array}
$$
By Bounded Covariant, for any $\tau\in\Omega_{k+7,m}^{I}$,
$$
\frac{|B_{\omega_0*\sigma_1*\tau_1*\tau}|}{|B_{\omega_0*\sigma_2*\tau_2*\tau}|}
\le\eta\frac{|B_{\omega_0*\sigma_1*\tau_1}|}{|B_{\omega_0*\sigma_2*\tau_2}|}
$$
By \eqref{gibcontract}, for $s=1,2$,
$$
1\le \frac{|B_{\omega_0}|}{|B_{\omega_0*\sigma_s*\tau_s}|}\le c_1^2.
$$
So, for any $\tau\in\Omega_{k+7,m}^{I}$,
$$
\frac{|B_{\omega_0*\sigma_1*\tau_1*\tau}|}{|B_{\omega_0*\sigma_2*\tau_2*\tau}|}
\le \eta c_1^2,
$$
which implies
$$\mu_{\beta,m}(B_{\omega_0*\sigma_1*\tau_1})\le \eta c_1^2
\mu_{\beta,m}(B_{\omega_0*\sigma_2*\tau_2}).$$ The case $\tau_1$
being of $D_2$,$\tau_2$ being of $\tilde{D}_2$ respectively(and
$\tau_1$ being of $D_3$,$\tau_2$ being of $\tilde{D}_3$ respectively
) can be discussed by the same way. Considering that, for $i=1,2,3$,
$$1\le\sharp D_i\le (2M+1)^3,\quad 1\le \sharp\tilde{D}_i\le (2M+1)^3,$$
we have
$$
\sum_{\tau\in D_i}\mu_{\beta,m}(B_{\omega_0*\sigma_1*\tau})\le (2M+1)^3 \eta c_1^2
\sum_{\tau\in \tilde{D}_i}\mu_{\beta,m}(B_{\omega_0*\sigma_2*\tau}).
$$
This implies that the inequality \eqref{progib4} holds.
\end{proof}

\begin{proof}[Proof of Theorem \ref{gibbs}]
We only prove the second inequality. Let $V$ be large enough so that Bounded covariation holds.

For any $k\ge1$, $\omega_0=i_0\cdots i_k\in \Omega_k$,
$m>k+3$, $\sigma=i_{k+1}i_{k+2}i_{k+3}\in \Omega_{k+1,k+3}^{i_k}$,
$$
\mu_{\beta,m}(\omega_0*\sigma)=b_m^{-1}\sum\limits_{\sigma_1\in\Omega_{k+4,m}^{i_{k+3}}}
|B_{\omega_0*\sigma*\sigma_1}|^\beta.
$$
So by \eqref{progib2}
$$
\mu_{\beta,m}(\omega_0)b_m\le
c\mu_{\beta,m}(\omega_0*\sigma)b_m
=c|B_{\omega_0*\sigma}|^\beta
\sum\limits_{\sigma_1\in\Omega_{k+4,m}^{i_{k+3}}}
\frac{|B_{\omega_0*\sigma*\sigma_1}|^\beta}{|B_{\omega_0*\sigma}|^\beta}.
$$
For any $\omega_1\in\Omega_{k+3}^{(k+3,i_{k+3})}$,
by bounded covariation,
$$%\begin{array}{rcl}
 \mu_{\beta,m}(\omega_0)b_m \le c
\eta^\beta|B_{\omega_0}|^\beta
\sum\limits_{\sigma_1\in\Omega_{k+4,m}^{i_{k+3}}}
\frac{|B_{\omega_1*\sigma_1}|^\beta}{|B_{\omega_1}|^\beta},$$ hence,
$$\mu_{\beta,m}(\omega_0) b_m |B_{\omega_1}|^\beta
\le c\eta^\beta|B_{\omega_0}|^\beta
\sum\limits_{\sigma_1\in\Omega_{k+4,m}^{i_{k+3}}}
|B_{\omega_1*\sigma_1}|^\beta.
%\end{array}
$$
Take sum on both side for any $\omega_1\in\Omega_{k+3}^{(k+3,i_{k+3})}$,
$$\mu_{\beta,m}(\omega_0) b_m
\sum\limits_{\omega_1\in\Omega_{k+3}^{(k+3,i_{k+3})}}|B_{\omega_1}|^\beta
\le c\eta^\beta|B_{\omega_0}|^\beta
\sum\limits_{\omega\in\Omega_{m}^{(k+3,i_{k+3})}}
|B_{\omega}|^\beta.$$ Take sum on both sides for all $i_{k+3}\in
\Omega_{k+3,k+3}$,
$$\mu_{\beta,m}(\omega_0) b_m b_{k+3} \le c\eta^\beta|B_{\omega_0}|^\beta b_m.$$
By \eqref{progib1},
$$\mu_{\beta,m}(\omega_0)\le c^2\eta^\beta b_{k}^{-1}|B_{\omega_0}|^\beta
= c^2\eta^\beta\mu_{\beta,k}(\omega_0).$$

Let $\mu_\beta$ be a weak limit of $(\mu_{\beta,m})_{m\ge1}$, we prove the theorem.
\end{proof}

\begin{proof}[Proof of Theorem \ref{haus}]
$(\mathscr{G}_n)_{n\ge0}$ is a sequence of coverings of $\sigma(H_{\alpha,V})$ with
diameter tends to $0$. So
$$\dim_H \sigma(H_{\alpha,V})\le s_*.$$
Now take any $\beta<s_*$, then $s_n>\beta$ for all $n$ large enough, thus
$$\sum_{B\in\mathscr{G}_n}|B|^{\beta}>
\sum_{B\in\mathscr{G}_n}|B|^{s_n}=1.$$
Let $\mu_\beta$ be a Gibbs-like measure defined in Theorem \ref{gibbs}.
Then for any large $k$ and each $B\in\mathscr{G}_k$ we have
$$\mu_\beta(B)\le\eta|B|^\beta.$$
Take $r>0$ small and $r$-Moran covering of $\sigma(H_{\alpha,V})$, i.e.,
$$\mathscr{M}=\{B\in\mathscr{G}_n\ :\ n\ge0, B\subset B^{-1}\in\mathscr{G}_{n-1},
|B^{-1}|>r, |B|\le r\}.$$
By Proposition A\ref{lm-2} and Theorem \ref{bvar}, for any $B\in\mathscr{M}$,
$$|B|>\frac{r}{\xi^2(V+5)(M+2)^3}.$$
For any ball $B(x,r)$, letting $\mathscr{C}=\{B\in\mathscr{M}: B\cap B(x,r)\ne\emptyset\}$,
$$\sharp \mathscr{C}\le 3\xi^2(V+5)(M+2)^3.$$
Then,
$$\mu_\beta(B(x,r))\le\sum_{B\in\mathscr{C}}\mu_\beta(B)
\le\eta\sum_{B\in\mathscr{C}}|B|^\beta
\le 3\eta \xi^2(V+5)(M+2)^3 r^\beta,
$$
which implies $\dim_H E>\beta$. Hence, $\dim_H E\ge s_*$.
\end{proof}

\bigskip

\noindent
{\bf Acknolegement}\ The authors thank Morningside Center of Mathematics
for the partial support.

\end{document}